\documentclass[11pt]{amsart}
\usepackage{amsthm,vmargin}
\usepackage{amsmath}
\begin{document}
\newcommand{\Q}{{\mathbb Q}}
\newcommand{\C}{{\mathbb C}}
\newcommand{\R}{{\mathbb R}}
\newcommand{\Z}{{\mathbb Z}}
\newcommand{\F}{{\mathbb F}}
\renewcommand{\wp}{{\mathfrak p}}
\renewcommand{\P}{{\mathbb P}}
\renewcommand{\O}{{\mathcal O}}
\newcommand{\Pic}{{\rm Pic\,}}
\newcommand{\Ext}{{\rm Ext}\,}
\newcommand{\rank}{{\rm rk}\,}
\newcommand{\sbull}{{\scriptstyle{\bullet}}}
\newcommand{\bX}{X_{\overline{k}}}
\newcommand{\ch}{\operatorname{CH}}
\newcommand{\tors}{\text{tors}}
\newcommand{\cris}{\text{cris}}
\newcommand{\alg}{\text{alg}}
\let\isom=\simeq
\let\rk=\rank
\let\tensor=\otimes

\newtheorem{theorem}[equation]{Theorem}      % (If you want theorem numbered
\newtheorem{lemma}[equation]{Lemma}          %
\newtheorem{corollary}[equation]{Corollary}  %       goes for lemmas, etc.)
\newtheorem{proposition}[equation]{Proposition}
\newtheorem{scholium}[equation]{Scholium}

\theoremstyle{definition}
\newtheorem{conj}[equation]{Conjecture}
\newtheorem*{example}{Example}
\newtheorem{question}[equation]{Question}

\theoremstyle{definition}
\newtheorem{remark}[equation]{Remark}

\numberwithin{equation}{section}

\title{Stability and locally exact differentials on a curve}
\author{Kirti Joshi}
\address{Math. department, University of Arizona, 617 N Santa Rita, Tucson
85721-0089, USA.}
\email{kirti@math.arizona.edu}
%%%%
%\date{Preliminary Version: May 26, 2003, revised: Feb 02, 2004, Feb 17, 2004}
%%%%%%%%%%%%%%%%%%%%%%%%%%%%%

\begin{abstract}
We show that the locally free sheaf $B_1\subset F_*(\Omega^1_X)$
of locally exact differentials on a smooth projective curve of
genus $g\geq 2$ over an algebraically closed field $k$ of
characteristic $p$ is a stable bundle. This answers a question of Raynaud.

\noindent\textsc{R\'ESUM\'E.}
%Nous montrons que le faiseau
%localement libre de differentielle localement exacte sur une
%courbe projective et lisse de genre $g> 1$ sur un corps
%alg\'ebriquement ferm\'e $k$ de  caract\'eristique $\neq 2$ est
%une fibre stable.
Soit $X$ une courbe propre,lisse, connexe, de genre $g$, d\'efinie
sur un corps $k$ alg\'ebriquement clos de caract\'eristique $p>0$.
Soit $F: X \to X$ le Frobenius absolu  et  $B_1\subset
F_*(\Omega^1_X)$, le faisceau des formes différentielles
localement exactes  sur X. C'est un fibr\'e vectoriel sur $X$ de
rang $p-1$. Nous montrons qu'il est stable pour $g\geq 2$.

\end{abstract}

\maketitle
\begin{verse}\textsf{\small
\hfill $\qquad\qquad\qquad$\b{dh}\=\i re \b{dh}\=\i re re man\=a,
 \b{dh}\=\i re sab ku\b{ch} hoye\\
\hfill $\qquad\qquad\qquad$ m\=ali
 si\d n\b{ch}e sau \b{gh}ar\=a, rit\=u \=aye \b{ph}al hoye\\
\hfill
$\qquad\qquad\qquad\qquad\qquad\qquad\qquad\qquad\qquad\qquad$
Kab\=\i r.}
\end{verse}

%\tableofcontents
\begin{center}

\end{center}

\section{Introduction}
          Let $k$ be an algebraically closed field of characteristic
$p>0$. Let $X/k$ be a smooth, projective curve of genus $g\geq 2$ over
$k$. Let $F:X\to X$ be the absolute Frobenius morphism of $X$.  If $V$
is a vector bundle we will write $\mu(V)=\deg(V)/\rk(V)$ for the slope
of $V$. We will say that $V$ is stable (resp. semi-stable) if for all
subbundles $W\subset V$ we have $\mu(W)<\mu(V)$ (resp. $\mu(W)\leq
\mu(V)$). We will write $\Omega^1_X=K_X$ for the canonical sheaf of
the curve.

Let $B_1$ be the locally free sheaf of locally exact differential
forms on $X$. This may also be defined by the exact sequence of
locally free sheaves
$$0\to \O_X\to F_*(\O_X) \to B_1 \to 0.$$ As $X$ is a curve, by
definition, $B_1$ is a sub-bundle of $F_*(\Omega^1_X)$ of rank $p-1$.
In \cite{raynaud82a} Raynaud showed that $B_1$ is a semi-stable bundle
of degree $(p-1)(g-1)$ and slope $g-1$. Raynaud has asked if $B_1$ is
in fact stable (see \cite{raynaud02a}). In this note we answer
Raynaud's question. We prove:
\begin{theorem}\label{main2}
Let $k$ be an algebraically closed field of characteristic $p>0$.
Let $X$ be a smooth, projective curve over $k$ of genus $g$ and
let $B_1$ be the locally free sheaf of locally exact differentials
on $X$. Then $B_1$ is a stable vector bundle of slope $g-1$ and
rank $p-1$.
\end{theorem}

Observe that when $p=2$, $B_1$ is a line bundle of degree $g-1$
and so the result is immediate in this case. Our proof also gives
a new proof of Raynaud's theorem that $B_1$ is semi-stable.
Theorem~\ref{main2} will follow from the corresponding assertion
for curves of sufficiently large genus. More precisely:

\begin{theorem}\label{main}
Let $k$ be an algebraically closed field of characteristic $p\neq
2$. Let $X$ be a smooth, projective curve over $k$ of genus $g$
and assume that the genus $g>(1/2)(p-1)(p-2)$. Then $B_1$ is
stable.
\end{theorem}

I thank M.~Raynaud for discussions which inspired me to think about
this problem again and for his comments. In the first version of this
note we had proved Theorem~\ref{main} (the large genus case).  That
this is sufficient to establish the result for all genus was pointed
out to us by Akio Tamagawa and I am grateful to him for his comments
and correspondence. This note was written while I was visiting Orsay
and I thank the mathematics department of Universit\'e de Paris-Sud
for support; thanks are also due to Jean-Marc Fontaine for many
stimulating conversations and support.

\section{The proofs}
We will first explain the reduction of Theorem~\ref{main2} to the
Theorem~\ref{main} (I owe this argument to A.~Tamagawa).

\begin{proof}{[Theorem~\ref{main} $\Longrightarrow$
Theorem~\ref{main2}]} As the genus $g_X$ of $X$ is at least two,
we know that there exists a connected finite \'etale covering of
$X$ of arbitrarily large genus.  Choose a finite \'etale covering
$f:Y\to X$ such that the genus $g_Y$ of $Y$  is sufficiently large
genus (more precisely with $g_Y>(p-1)(p-2)/2$). Observe that the
formation of $B_1$ commutes with any finite \'etale base change,
that is, $B_{1,Y}=f^*(B_{1,X})$. Hence the stability of $B_{1,X}$
follows from that of $B_{1,Y}$ and the latter assertion is the
content of Theorem~\ref{main}.
\end{proof}

The rest of this note will be devoted to the proof of
Theorem~\ref{main}. The idea of the proof is to get an upper bound
on the slope of the destabilizing subsheaf (if it exists!). We
recall some facts which we need. The first two lemmas are from
\cite{joshi02a} which is not yet published so we provide proofs.

\begin{lemma}\label{degree-formula}
We have $\deg(F_*(\Omega^1_X))=(p+1)(g-1)$.
\end{lemma}
\begin{proof} By Riemann-Roch theorem
$\chi(\Omega^1_X)=\chi(F_*(\Omega^1_X))$. Which gives
$\deg(F_*(\Omega^1_X))+p(1-g)=\deg(\Omega^1_X)+(1-g)$. Simplifying
this gives $\deg(F_*(\Omega^1_X))=2(g-1)+p(g-1)+(1-g)$, which
leads to the claimed result.
\end{proof}

\begin{lemma}\label{line-estimate}
Let $M\subset F_*(\Omega^1_X)$ be any line subbundle. Then
$$\mu(M)\leq
\mu(F_*(\Omega^1_X))-\frac{(p-1)(g-1)}{p}$$
\end{lemma}

\begin{proof}
This was proved in \cite{joshi02a}. By adjunction we get a map
$$F^*(M)\to \Omega^1_X.$$ Hence by degree considerations we see that
$$\mu(M)\leq
\frac{\deg(\Omega^1_X)}{p}=\mu(F_*(\Omega^1_X))-\frac{(p-1)(g-1)}{p},$$
where we use Lemma~\ref{degree-formula} for the last equality.
\end{proof}

We will also need to identify the dual $B^*_1$ of $B_1$. This was done
by Raynaud in \cite{raynaud82a}.

\begin{lemma}
The dual $B_1^*=B_1\tensor K_X^{-1}$.
\end{lemma}

Now we recall a theorem of Mukai-Sakai (see \cite[page
251]{mukai85a}). This will be used to give upper bounds on
destabilizing subbundles along with Lemma~\ref{line-estimate}. For
more on this see Remark~\ref{a remark}.

\begin{theorem}\label{mukai-thm}
Let $X/K$ be a smooth, projective curve over an algebraically closed
field $K$ of arbitrary characteristic. Let $W$ be a vector bundle. Fix
an integer $1\leq k\leq r=\rk(W)$. Then there exists a subbundle
$U\subset W$ of rank $k$ such that
$$\mu(W)\leq \mu(U)+g(1-k/r).$$
\end{theorem}

\begin{proof}{[of Theorem~\ref{main}]}
As remarked in the introduction, when $p=2$, $B_1$ is a line
bundle and so there is nothing to prove in this case. In what
follows we will assume that $p\neq2$. Assume, if possible, that
$B_1$ is not stable. Then there exists some $W\subset B_1$ which
is locally free of rank $r$ with $\mu(W)\geq \mu(B_1)=g-1$.

We first claim that we can assume with out loss of generality that
$r\leq (p-1)/2$. Indeed, if not then the dual of $B_1$ surjects on
$W^*$ the dual of $W$ and the kernel has rank $\leq (p-1)/2$. Moreover
writing $W_1$ for the kernel we get an exact sequence
$$0\to W_1\to B_1^*=B_1\tensor K_X^{-1} \to W^* \to 0$$
and writing $W_2=W_1\tensor K_X$, we get
$$0\to W_2\to B_1\to W^*\tensor K_X\to 0$$
Now $\deg(W_1)+\deg(W^*)=\deg(B_1^*)=(p-1)(1-g)$ and a simple
calculation using $-\mu(W^*)=\mu(W)\geq (g-1)$ shows that
$\deg(W_1)\geq \rk(W_1)(1-g)$.
%so that
%\begin{eqnarray}
%\deg(W_1)&=&(p-1)(1-g)-\deg(W^*)\\
%         &\geq&(p-1)(1-g)-(1-g)\rk(W^*)\\
%         &\geq&[(p-1)-\rk(W^*)](1-g)\\
%         &\geq&\rk(W_1)(1-g),
%\end{eqnarray}
% where we have used $-\mu(W^*)=\mu(W)\geq (g-1)$
So that we have $$\deg(W_2)=\deg(W_1\tensor
K_X)\geq\rk(W_1)(1-g)+\rk(W_1)(2g-2),$$ which simplifies to
$$\deg(W_2)\geq\rk(W_2)(g-1)$$ or equivalently $\mu(W_2)\geq(g-1)$.

Thus we may assume without loss of generality that $\rk(W)\leq
(p-1)/2$.

We apply Theorem~\ref{mukai-thm} to the following situation. We take
$W$ as above of slope $\geq(g-1)$ with $\rk(W)$ at most $(p-1)/2$ and
we take $k=1$ and let $U$ to be the line bundle given by
Theorem~\ref{mukai-thm}. Then we get
$$\mu(W)\leq \mu(U)+g(1-1/\rk(W)).$$ Now as $U\subset W\subset
B^1\subset F_*(\Omega^1_X)$, we know by Lemma~\ref{line-estimate} that
$$\mu(U)\leq \mu(F_*(\Omega^1_X))-(p-1)(g-1)/p.$$

Putting these two inequalities together and using $\mu(W)\geq g-1$
we get
\begin{equation}
(g-1)\leq\mu(W)\leq \mu(F_*(\Omega^1_X))-(p-1)(g-1)/p+g(1-1/\rk(W)).
\end{equation}
Writing $r=\rk(W)$ this simplifies to
\begin{eqnarray}
(g-1)\leq\mu(W)&\leq& \frac{2(g-1)}{p}+g(1-1/r)\\
&\leq& \left(\frac{2}{p}+1-\frac{1}{r}\right)g-\frac{2}{p}
\end{eqnarray}
And so in particular we get
$$g-1\leq \left(\frac{2}{p}+1-\frac{1}{r}\right)g-\frac{2}{p},$$
which easily simplifies to
%Collecting all the terms involving $g$ together we get
%$$g\left(\frac{1}{r}-\frac{2}{p}\right)\leq 1-\frac{2}{p}$$
%or that
%$$
%g<\left(\frac{p-2}{p}\right)/\left(\frac{p-2r}{pr}\right)
%=\left(\frac{p-2}{p}\right)\left(\frac{pr}{p-2r}\right).$$
%or
$$g \leq \frac{r(p-2)}{p-2r},$$ and as $r$ can be as large as $(p-1)/2$,
this says that $g\leq\frac{1}{2}(p-1)(p-2)$.  Now we have assumed that
$g>\frac{1}{2}(p-1)(p-2)$ so that our assumption that $W\subset B_1$
has slope $\geq g-1$ leads to a contradiction. Thus $B_1$ is stable.
\end{proof}

\begin{remark}
In fact, as was pointed out to us by A.~Tamagawa,
Theorem~\ref{main} can be slightly sharpened as follows. We claim
that  the following conditions are equivalent
\begin{enumerate}
\item The bundle $B_1$ is stable.
\item Either $p=2$ or $g>1$.
\end{enumerate}
After Theorem~\ref{main2}, the assertion (2) $\Longrightarrow$ (1)
is immediate. Now to prove that (1) $\Longrightarrow$ (2) observe
that the case $p=2$ is trivial so assume, if possible, that $g\leq
1$. By Theorem~\ref{mukai-thm} for line bundles applied to $B_1$,
we see that $B_1$ has a line subbundle $M\hookrightarrow B_1$ with
$$\mu(M) \geq \mu(B_1)-g\left(1-\frac{1}{p-1}\right).$$  Then the right hand side is greater than
$(g-1)-1$. As $\mu(M)$ is an integer, this says that $\mu(M)\geq
(g-1)$ and this contradicts the  stability of $B_1$.
\end{remark}

\begin{remark}\label{a remark}
The method of using coverings of large degree can also be used to
improve the results of \cite{joshi03a} (resp. \cite{joshi02a} for
small $p$) where we showed that $F_*(L)$ is stable for any line
bundle $L$ on $X$ provided that the genus of $X$ is sufficiently
large.
\end{remark}
\begin{remark}
The method of proof given above can also be used to show that the
bundles $B_{n}$ defined in \cite{illusie79b} using iterated
Cartier operations are stable. These are bundles of rank $p^n-1$
and slope $g-1$.
\end{remark}

%\bibliographystyle{amsalpha}
%\bibliography{raynaud}

\begin{thebibliography}{JRXY02}

\bibitem[Ill79]{illusie79b}
L.~Illusie, \emph{Complexe de de {R}ham-{W}itt et cohomologie cristalline},
  Ann. Scient. Ecole Norm. Sup. \textbf{12} (1979), 501--661.

\bibitem[Jos03]{joshi03a}
K.~Joshi, \emph{On vector bundles destabilized by frobenius {III}}, Preprint
  under preparation (2003).

\bibitem[JRXY02]{joshi02a}
K.~Joshi, {S}. {R}amanan, {E}.~{Z}. {X}ia, and {J}.~{K}. Yu, \emph{On vector
  bundles destabilized by frobenius}, Preprint (2002).

\bibitem[MS85]{mukai85a}
S.~Mukai and {F}. Sakai, \emph{Maximal subbundles of vector bundles on a
  curve}, Manuscripta {M}ath. \textbf{52} (1985), 251--256.

\bibitem[Ray82]{raynaud82a}
M.~Raynaud, \emph{Sections des fibr\`es vectoriels sur une courbe}, Bull.
  {S}oc. {M}ath. {F}rance \textbf{110} (1982), no.~1, 103--125.

\bibitem[Ray02]{raynaud02a}
Michel Raynaud, \emph{Sur le groupe fondamental d'une courbe compl\`ete en
  caract\'eristique {$p>0$}}, Arithmetic fundamental groups and noncommutative
  algebra (Berkeley, CA, 1999), Proc. Sympos. Pure Math., vol.~70, Amer. Math.
  Soc., Providence, RI, 2002, pp.~335--351.

\end{thebibliography}
\providecommand{\bysame}{\leavevmode\hbox to3em{\hrulefill}\thinspace}
\providecommand{\MR}{\relax\ifhmode\unskip\space\fi MR }
\providecommand{\MRhref}[2]{%
  \href{http://www.ams.org/mathscinet-getitem?mr=#1}{#2}
}
\providecommand{\href}[2]{#2}

\end{document}